\documentclass [11pt]{article}
\usepackage{amsmath, amsthm, amssymb}
\usepackage{graphicx}
\newfont{\bbb} {msbm10}
\newcommand{\R}{\Bbb{R}}
\newcommand{\bS}{\Bbb{S}}
\newcommand{\bB}{\Bbb{B}}
\newcommand{\T}{\Bbb{T}}

\newcommand{\cT}{{\cal{T}}}

\newcommand{\sbs}{\subset}
\newcommand{\ra}{\rightarrow}

\newcommand{\cP}{{\cal{P}}}

\newcommand{\p}{\partial}

\newcommand{\pt}{\frac{\partial}{\partial t}}
\newcommand{\n}{\nabla}

\newcommand{\cG}{{\cal{G}}}
\newcommand{\cA}{{\cal{A}}}

\newcommand{\B}{\Bbb{B}}

\newcommand{\HH}{\Bbb{H}}

\newcommand{\0}[1]{_{_{#1}}}

\usepackage[all]{xy}

\begin{document}

\title{On the Farrell and Jones Warping Deformation}
\author{Pedro Ontaneda\thanks{The author was
partially supported by a NSF grant.}}
\date{}

\maketitle

\begin{abstract} 
The Farrell-Jones warping deformation is a powerful geometric construction that has been crucial
in the proofs of many important contributions to the theory of manifolds of negative
curvature. In this paper we study this construction in depth, in a more general setting, and obtain 
explicit quantitative results.

The results in this paper are key ingredients in the problem of smoothing Charney-Davis strict
hyperbolizations  \cite{ChD}, \cite{O1}.\vspace{.1in}

\end{abstract}

\noindent {\bf \large  Section 0. Introduction.}

Let $g$ be a Riemannian metric on the $n$-sphere
$\bS^n$. Consider the warped metric $h=sinh^2(t)\, g +dt^2$ on $\R^{n+1}-\{0\}=\bS^n\times (0,\infty)$. If $g=\sigma\0{\bS^n}$, the canonical round metric on $\bS^n$, then $h$ is the (real) hyperbolic metric.  
But for general $g$ the metric $h$ is not hyperbolic. In \cite{FJ1} Farrell and Jones used the following
method to deform the metric $h$ to a hyperbolic metric in a ball of large radius $2\alpha$ centered at the origin.\vspace{.1in}

For $\alpha>0$ consider the metric: $$h\0{\alpha}(x,t)=
sinh^2(t)\,\Big( \big(1-\rho\0{\alpha}(t)\big)\, \sigma\0{\bS^n}(x)\,\,+\,\, \rho\0{\alpha}(t)\,g(x) \Big)\,\,\,+\,\,\,dt^2$$
\noindent where $\rho\0{\alpha}(t)=\rho(\frac{t}{\alpha}-1)$, and $\rho:\R\ra [0,1]$ is
a fixed smooth function with $\rho(t)=0$ for $t\leq 0$ and $\rho(t)=1$ for $t\geq 1$.
Hence, for $t\leq \alpha$, the metric $h\0{\alpha}$ is hyperbolic, for $t\geq 2\alpha$
we have $h\0{\alpha}=h$ and between $t=\alpha$ and $t=2\alpha$ the metric
$\sigma\0{\bS^n}$ deforms to $g$. The metrics $h\0{\alpha}$ have two important
properties: 

\begin{enumerate}
\item[(1)] they are all hyperbolic for $t\leq\alpha$, i.e. in the ball of radius $\alpha$ centered at the origin,
\item[(2)] given $\epsilon>0$ there is $\alpha\0{0}$ such that all sectional curvatures of 
$h\0{\alpha}$ lie within $\epsilon$ of -1, provided 
$\alpha>\alpha\0{0}$. This was proved in \cite{FJ1}.\end{enumerate}

\noindent {\bf Remark 0.1.} Roughly speaking, the reason (2) holds is because taking $\alpha$ large has two consequences:
(a) the deformation from $\sigma\0{\bS^n}$ to $g$ happens far from the origin, i.e. when $t$ is large: $t>\alpha$, and (b) if $\alpha$ is large 
the deformation from $\sigma\0{\bS^n}$ to $g$ happens slowly, i.e when $t\in [\alpha,2\alpha]$.\vspace{.1in}

As mentioned in the abstract this warping deformation  has been crucial
in the proofs of many important contributions to the theory of manifolds of negative
curvature (see for instance  \cite{AF1}, \cite{AF2}, \cite{AF3}, \cite{A}, \cite{FJ1}, \cite{FJ1.5}, \cite{FJ2},
 \cite{FJO1},  \cite{FJO2},
\cite{FOR}, \cite{FO1}, \cite{FO1.5}, \cite{FO2}, \cite{FO3}, \cite{FO4}, \cite{FO5},  \cite{O}).
To state our results we need some concepts and definitions.
\vspace{.1in}

Let $g\0{t}$ be a smooth one-parameter family of metrics on a closed manifold $M^n$, and consider the
metric  $h=sinh^2(t)\, g\0{t} +dt^2$ on $M\times I$, where $I\sbs (0,\infty)$ is an interval. We say that 
$\{ g\0{t}\}$ is {\it  $\epsilon$-slow} if $g\0{t}$ and its first derivatives in the $M$-direction
change $\epsilon$-slowly with $t$ (up to order 2 and 1, respectively; 
see Section 3 for an intrinsic definition). 
A metric $g$ on a compact manifold is {\it c-bounded} if the derivatives up to
order 2 of $g$ are bounded by $c$ (and $|det \,g|>1/c$, see Section 3 for more details). A set of metrics $\{ g\0{t}\}$  is {\it c-bounded} if every $g\0{t}$
is $c$-bounded.   Our Main Theorem says that
if $\{g\0{t}\}$ is $\epsilon$-slow and $c$-bounded then
the metric $h$ is $\eta$-close to being hyperbolic, where $\eta=\eta(\epsilon, c, n)$.
The concept of a metric being {\it $\eta$-close to hyperbolic} is given in the next paragraph
(for more details see Section 2).\vspace{.1in}

Let $\bB$ be the unit $n$-ball, 
with the flat metric $\sigma\0{\R^n}$.  
Our basic model is $\T=\bB\times (-1,1)$,  with hyperbolic metric $\sigma=e^{2t}\sigma\0{\R^n}+dt^2$. Let $(M,g)$ be a Riemannian manifold and $S\sbs M$. We say that $S$ is $\epsilon$-{\it close to hyperbolic} if for every $p\in S$ there is an {\it $\epsilon$-close to hyperbolic 
chart with center $p$}, that is, there is a chart
$\phi :\T\ra M$, $\phi(0,0)=p$,  such that $|\phi^*g-\sigma|\0{C^2}<\epsilon$,
where $|.|\0{C^2}$ is the $C^2$ norm.  \vspace{.1in}

If $M=N\times J$, $J$ an interval,  we say that $S\sbs M$ is {\it radially $\epsilon$-close to hyperbolic} if, in addition, for every $p\in S$ there is an  $\epsilon$-close to hyperbolic 
chart $\phi$ with center $p$ and, in addition, the chart 
$\phi$ respects the product structure of $\T_\xi$
and $M$,
i.e. $\phi(. , t)=(\phi_1(.), t+a)$, for some $a$ depending on the
$\phi$ (see Section 2  for details). Here the ``radial" directions
are $(-1,1)$ and $J$ in $\T$ and $M$, respectively.
\vspace{.1in}

\noindent{\bf Remark 0.2.} 
The definition of radially $\epsilon$-close to hyperbolic metrics 
is well suited to studying warp metrics for $t$ large, but for
small $t$ this definition is not useful because of: (1) the need for some space to fit the charts, and (2) the form of our specific fixed model $\T$.
An undesired consequence is that punctured hyperbolic space $\HH^n-\{ o\}=\bS^{n-1}\times\R^+$
(with warp metric $sinh^2(t)\sigma\0{\bS^{n-1}}+dt^2$) is not radially $\epsilon$-close to hyperbolic
for $t$ small. In fact there is ${\sf a}={\sf a}(n, \epsilon)$ such that hyperbolic
space is $\epsilon$-close to hyperbolic for $t>{\sf a}$ (and not for all $t\leq {\sf a}$). We prove this in Section 3, where we also give an explicit formula for
${\sf a}$ (see 3.9). 
\vspace{.1in}

Now we can state our main result.\vspace{.1in}

\noindent {\bf Main Theorem.}  {\it Let $M^n$ be a closed smooth manifold
and $I=(a,b)\sbs (0,\infty)$ an interval.
If the family of metrics $\{g\0{t}\}_{t\in I}$ on $M$ is $\epsilon$-slow and $c$-bounded then
the metric  $h=sinh^2t\,g\0{t}+dt^2$ is radially $\eta$-close to hyperbolic on $M\times I'$,
provided $$\,C_1\,(e^{-a}+\epsilon)\leq\eta$$ 
\noindent where $C_1=C_1(c,n)$ and $I'=(a+1, b-1)$.}\vspace{.1in}

\noindent {\bf Remarks 0.3.} 

\noindent {\bf 1.} An explicit formula for $C_1$ is given at the end of the proof of Corollary 3.3.

\noindent {\bf 2.} We actually prove a slightly more general result 
(Corollary 3.3) in which the ``size of the charts" 
varies (depends on a variable $\xi$). In this case the constant $C_1$ also depends on $\xi$.

\noindent {\bf 3.} The construction of the $\eta$-close to hyperbolic
charts around a given point is explicit; see formula $(\star)$ in the
proof of Theorem 3.2. 
\vspace{.1in}

It is interesting to compare the statement in the Main Theorem with the Farrell-Jones deformation
mentioned at the beginning of the Introduction (see also Remark 0.1).
Write $g\0{t}=(g\0{\alpha})\0{t}= (1-\rho\0{\alpha}(t))\, \sigma\0{\bS^n}(x)\,\,+\,\, \rho\0{\alpha}(t)\,g(x)$. Then if $\alpha$ is large the family $g\0{t}$ varies slowly. Also, since $\bS^n$ is
compact the family $g\0{t}$ is bounded. Note that the Main Theorem says more than that 
$h=sinh^2t\,(g\0{\alpha})\0{t}+dt^2$ has curvatures close to -1, it says that (in a chart sense) the metric $h$
is close to being hyperbolic. And, more importantly, we get an explicit relationship between
how large $t$ has to be (see (a) in Remark 0.1), how slowly $g\0{t}$ varies (see (b) in Remark 0.1)
and how close the metric $h$ gets to being hyperbolic. Also note that in the Main Theorem the result
holds for $t>a+1$, so the variable $a$ tells us how large $t$ has to be. Therefore, if $t$ is very large 
and $g\0{t}$ varies very slowly, then $h$ is very close to being hyperbolic; and how close $h$ gets to being
hyperbolic is given
by the formula in the Main Theorem: it decreases exponentially in terms of how large $t$ is, and linearly in
terms of how slowly $g\0{t}$ varies.\vspace{.1in}

Let us consider again the setting of the Farrell-Jones deformation. 
As before let $g$ be a metric on the $n$-sphere $\bS^n$ and $\rho$ as above.
Given positive numbers $a$ and $d$ define $\rho\0{a,d}(t)=\rho(2\,\frac{t-a}{d})$.
Note that $\rho\0{\alpha,2\alpha}=\rho\0{\alpha}$.
Write $(g\0{a,d})\0{t}=\sigma\0{\bS^n} +\rho\0{a,d}(t) (g-\sigma\0{\bS^n})$ and define the metric
\begin{center} $
\cT_{_{a,d}}\, g\, =\, sinh^2\, t\,\,(g\0{a,d})\0{t}+dt^2
$\end{center}
\noindent Sometimes we will also write
$ \cT_{_{a,d}}\, h$ instead of $\cT_{_{a,d}}\, g$, where
$h=sinh^2(t) g+dt^2$. By construction we have\vspace{.12in}

\noindent (*)\hspace{1.3in}$\cT\0{a,d}g\,=\,\left\{ \begin{array}{lllll}
sinh^2\, (t)\sigma\0{\bS^{n}}\, +\, dt^2&& {\mbox{on}}& & B_{a}\\
g&&{\mbox{outside}}& & B_{a+\frac{d}{2}}
\end{array}\right.$\vspace{.12in}

\noindent  where $B_a\sbs\R^{n+1}$ is the ball of radius $a$ centered at the origin.
Hence inside the ball $B_a$ the metric $\cT\0{a,d}g$ is hyperbolic.
We call the process $g\mapsto \cT_{_{a,d}}\, g$ 
the {\it two variable warping deformation}. (The two variables are $a$ and $d$). 
Note that  $\cT_{_{\alpha,2\alpha}}\, g$ coincides with the metric $h\0{\alpha}$ given by the Farrell-Jones deformation, mentioned
at the beginning of this introduction. Also note that one of  the differences
with the Farrell-Jones deformation is that now ``how large $t$ has to be" (given by $\alpha$ in the Farrell-Jones deformation
and by $a$ in $\cT_{_{a,d}}\, g$) and ``how long we have to stretch the change from $\sigma\0{\bS^n}$ to $g$"
(given also by $\alpha$ in the Farrell-Jones deformation and by $d/2$ in $\cT_{_{a,d}}\, g$) are independent
variables. The following result will be deduced from the Main Theorem in Section 4. As before let $a,\, d>0$; also $b>1$.
\vspace{.1in}

\noindent {\bf Theorem 1.} {\it Let the metric $g$ on $\bS^n$ be $c$-bounded.  Then \vspace{.1in}

\noindent (1) the metric $\cT_{_{a,d}}\, g\,$ is 
hyperbolic on $B_{a}$ 

\noindent (2) the metric $\cT_{_{a,d}}\, g\,$ is 
radially $\epsilon$-close to hyperbolic outside $B_{b}$,\vspace{.1in}
 
\noindent provided} 
\begin{center}$\,C_2\,\Big(e^{-b}+\frac{1}{d}\,\Big)\leq\epsilon$
\end{center}
\noindent {\it where $C_2=C_2(c,n)$. } \vspace{.1in}

\noindent {\bf Remarks.}

\noindent{\bf 1.} 
A formula for $C_2(c,n)$ is given in Section 4. 

\noindent{\bf 2.} The term``radially" in Statement (2) refers to
the decomposition $\R^{n+1}-\{ 0\}=\bS^n\times(0,\infty)$.

\noindent{\bf 3.}  
In Theorem 1 we do not say that $\cT_{_{a,d}} g$
is $\epsilon$-close to hyperbolic inside $B_{b}$. 
This is because $\cT_{_{a,d}} g$ is not $\epsilon$-close to hyperbolic for $t$ small  (see Remark 0.2).

\noindent{\bf 4.} As with the Main Theorem we actually prove a slightly  more general result in which the ``size of the charts " depends on a variable $\xi$. In this case the constant
$C_2$ also depends on $\xi$. \vspace{.1in}

Theorem 1 is most useful when
$b<a$. In this case for any $p\in\R^{n+1}$ the metric $\cT_{_{a,d}} g$ is either hyperbolic near $p$ or it is radially $\epsilon$-close to hyperbolic
near $p$. This motivates the following definition. 
Let $B_{a}=B_a(0)$ be the ball of radius $a$ centered at $0$. We say
that a metric $h$ on $\R^{n+1}$ is $(B_a,\epsilon)$-{\it close to
hyperbolic} \, if \vspace{.1in}

\begin{enumerate}
\item[  (1)]  On $B_{a}-\{0\}=\bS^n\times (0,a)$
we have $h=sinh^2(t)\sigma\0{\bS^n}+dt^2$. Hence $h$
is hyperbolic on $B_a$.
\item[(2)]  the metric $h$ is 
radially $\epsilon$-close to hyperbolic outside $B_{a-1}$.
\end{enumerate}\vspace{.1in}

\noindent {\bf Remarks.}

\noindent {\bf 1.} We have dropped the word ``radially" 
to simplify the notation. But it does appear in condition (2),
where now ``radially" refers to the center on $B_a$.

\noindent{\bf 2.} We will always assume $a>{\sf a}+1$,
where {\sf a} is as in 0.2. Therefore conditions (1), (2) 
and remark 0.2 imply
a stronger version of (2):\vspace{.1in}

  (2') the metric $h$ is 
radially $\epsilon$-close to hyperbolic outside $B_{\sf a}$.\vspace{.1in}

\noindent This is the reason why we demanded radius $a-1$
in (2), instead of just $a$. 

\noindent {\bf 3.} In the slightly more general case where
the ``size $\xi$ of charts varies we have to take
$a>{\sf a}+1+\xi$. \vspace{.1in}

Metrics that are $(B_a,\epsilon)$-close to
hyperbolic are very useful, and are key objects in
\cite{O1}. See also \cite{O2}, \cite{O4}.
Hence it is helpful to have some notation for this type of metrics.
For instance, with this new notation, Theorem
1 can be restated in the following way (taking
$b=a-1$): \vspace{.1in}

\noindent {\bf Corollary.} {\it Let the metric $g$ on $\bS^n$ be $c$-bounded. Then the metric $\cT_{_{a,d}}\, g\,$ is  $(B_a,\epsilon)$-close to hyperbolic provided} 
\begin{center}$\,C'_2\,\Big(e^{-a}+\frac{1}{d}\,\Big)\leq\epsilon$
\end{center}
\noindent {\it where $C'_2=e\,C_2$, $C_2$ as in Theorem 1. } \vspace{.1in}

Now let us fix $\epsilon$ first and take: (1) the number $a$ such
that $C_2 e^{-a}<\epsilon/2$ (and $a>{\sf a}$) and (2) the number
$d$ such that $C_2/d<\epsilon/2$. Applying these choices to
Corollary we obtain: \vspace{.1in}

\noindent {\bf Theorem 2.} {\it Let the metric $g$ on $\bS^n$ be $c$-bounded and $ \epsilon>0$. Then the metric $\cT_{_{a,d}}\, g\,$ is 
$(B_a,\epsilon)$-close to hyperbolic provided we take
$a$ and $d$ large enough. Explicitly, we have to take} 
\begin{center}$a\,>\, a(c,\epsilon,n)\hspace{.4in} and\hspace{.4in}  d\,>\, d
(c,\epsilon,n)$\end{center}\vspace{.1in}

\noindent {\bf Remarks.}

\noindent {\bf 1.} We can take  $a(c,\epsilon,n)=ln(\frac{2C_2}{\epsilon})+{\sf a}(\epsilon, n+1)$
and $d(c,\epsilon,n)=\frac{2C_2}{\epsilon}$. 

\noindent{\bf 2.}  We actually prove a slightly more general result in which the ``size of the charts" depends on a variable $\xi$. In this case the constants $a(c,\epsilon,n)$ and $d(c,\epsilon,n)$ also depend on $\xi$. \vspace{.1in}

\vspace{.1in}

The results in this paper are key ingredients in the problem of smoothing Charney-Davis strict
hyperbolizations  \cite{ChD}, \cite{O1}. Next we give an idea how the two variable warping deformation fits in the smoothing problem.\vspace{.1in}



In the same way that a cubical complex is made of basic pieces (the cubes
$\square^k$),
the hyperbolization $h(K)$ of a cubical complex $K$ is 
also made of
basic pieces: pre-fixed hyperbolization pieces $X^k$. Indeed one begins with a cubical complex $K$ and replaces each cube of dimension $k$
by the hyperbolization piece of the same dimension. Cube complexes have
a piecewise flat metric induced from the flat geometry of the cubes.
Likewise the Charney-Davis hyperbolizations have a piecewise hyperbolic metric
 because the Charney-Davis hyperbolization pieces
are hyperbolic manifolds (compact, with boundary and corners).
To see how singularities appear one can first think about the manifold 2-dimensional 
cube case. If $K^2$ is a 2-dimensional manifold cube complex then
its piecewise flat metric is Riemannian outside the vertices. A vertex is
a singularity if and only if the vertex does not meet exactly four cubes.
The picture is exactly the same for $h(K^2)$. 
These point singularities in $h(K^2)$ can be smoothed out
easily using warping methods.
In higher dimensions
the singularities of $K^n$ and $h(K)$ appear in (possibly the whole of)
the codimension 2 skeletons $K^{(n-2)}$ and $h(K^{(n-2)})$, respectively.
In \cite{O1} the idea of smoothing the piecewise hyperbolic metric on
$h(K)$ is to do it inductively down the dimension of the skeleta.
One begins with the $(n-2)$-dimensional pieces $X^{n-2}$. Transversally
to each $X^{n-2}$ (that is, on the union of geodesic segments emanating
perpendicularly to $X^{n-2}$, from a fixed point in $X^{n-2}$) one has
essentially the 2-dimensional picture mentioned above. Once we
solve this transversal problem we extend this transversal smoothing by taking a warp product
with $X^{n-2}$; we called this product method {\it hyperbolic extension}
\cite{O2}.
This gives a smoothing on a (tubular) neighborhood of the piece $X^{n-2}$.
Caveat: we do not want to actually have a smoothing on a neighborhood
of the whole of $X^{n-2}$, since we will certainly have 
matching problems for different $X^{n-2}$ meeting on a common
$X^{n-3}$; so we only want a smoothing on a neighborhood the the $Z^{n-2}$,
where $Z^{n-2}\sbs X^{n-2}$ is just a bit ``smaller" than $X^{n-2}$,
so that the neighborhoods of the $Z^{n-2}$ are all disjoint. Next 
step is to smooth around the $X^{n-3}$ (or, specifically the $Z^{n-3}$).
The metric is already smooth outside a neighborhood of the $(n-3)$-skeleton. Transversally to each $X^{n-3}$ we
have a 3 dimensional problem.
(It helps to have a 3 dimensional picture in mind, like in dimension 2).
It happens that if we did things with care in the first step (around the
$Z^{n-2}$) the metric  in the 3 dimensional transversal problem
(let's call this metric $\cP\0{3}$) is radially $\epsilon$-close to hyperbolic outside some large ball B.
Here is where we want to use two variable warping deformation
given in this paper:
to extend $\cP\0{3}$ to  a metric $\cG\0{3}$ defined on the whole of B. But we have
a problem because $\cP\0{3}$ is radially $\epsilon$-close to hyperbolic
but not warped. To solve this problem we apply {\it warp forcing} \cite{O3}
to $\cP\0{3}$ first, to obtain a warp metric $sinh^2(t)g\0{3}+dt^2$.
This warp metric is still close to hyperbolic (with a larger but
controllable $\epsilon$, see \cite{O3}). Now we apply the two variable
warping deformation to the warp metric $sinh^2(t)g\0{3}+dt^2$
to obtain our desired extension $\cG\0{3}$ of $\cP\0{3}$, getting rid, in this way, of the transverse
singularity. Note that $\cG\0{3}$ is also 
close to hyperbolic (with a larger but
controllable $\epsilon$, see Theorem 2 above). Once the transversal 3 dimensional problem
is solved we extend this smoothing to neighborhoods of the $Z^{n-3}$
using hyperbolic extension. Next we do the same for the $Z^{n-4}$
and so on until we smooth out all the singularities.\vspace{.1in}

The paper has 4 sections. In Section 1 we give some notation. In Section 2 we introduce the standard
models and $\epsilon$-close to hyperbolic metrics. In Section 3 we define $\epsilon$-slow families of
metrics and prove (the slightly more general version of) the Main Theorem. In Section 4 we deal with the
two variable warping deformation and prove Theorem 1.

\vspace{.2in}

\noindent {\bf \large Section 1.  Notation.}

Let $f:M_2\ra \R^+=(0,\infty)$ be smooth. Recall that the metric $g=f^2g_1+g_2$ is called a {\it warped metric} on $M_1\times M_2$, and $f$ is the 
{\it warping function}. 
(For a study of warped metrics see \cite{BisOn}.) In the one-dimensional
case, that is, when $M_2=I\sbs\R$ is an interval, the warped metric on $M\times I$ is written $f^2g+dt^2$, where
$f:I\ra \R^+$ and $g$ is a metric on $M$.\vspace{.1in}

Let $g_t$, $t\in I\sbs\R$, be a one-parameter family of metrics on the manifold $M$.
(All one-parameter families in this paper will be assumed to be smooth, were ``smooth" means $(x,t)\mapsto g(x,t)=g_t\mid_x$ is smooth.)
 We call the metric $h=g_t+dt^2$ on $M\times I$ a {\it variable metric
with metrics $g_t$}. Note that the $t$-lines $t\mapsto (x,t)$, $x\in M$, are the integral curves of the vector field $\pt$ on $M\times I$. It follows from Koszul's formula that the t-lines are geodesics. \vspace{.1in}

On a Riemannian manifold $(N^n,h)$ we
can identify, using the exponential map, a geodesic ball $B(p,\epsilon)$, with its center $p$  deleted, with the cylinder $\bS^{n-1}\times (0,\epsilon)$.
Then, by Gauss Lemma, $h$ is a variable metric on this cylinder, and we write $h=g_r+dr^2$, where $r$ is the distance to $p$.\vspace{.1in}

Note that every warped metric $f^2g+dt^2$ on $M\times I$ is a variable metric on $M\times I$, with $g_t=f^2(t)g$.


\vspace{.2in}

\noindent {\bf \large Section 2. The basic local hyperbolic model and $\epsilon$-close to hyperbolic metrics.}

Let $\bB=\bB^n\sbs\R^n $ be the unit ball, 
with the flat metric $\sigma\0{\R^n}$.  Write  $I_\xi=(-(1+\xi),1+\xi)\sbs\R$, $\xi\geq0$.
Our basic models are $\T^{n+1}_\xi=\T_\xi=\bB\times I_\xi$,  with hyperbolic metric $\sigma=e^{2t}\sigma\0{\R^n}+dt^2$. 
In what follows we may sometimes suppress the sub index $\xi$, if the context is clear.
The number $\xi$ is called the {\it excess} of $\T_\xi$.
\vspace{.1in}

\noindent {\bf Remarks.}

\noindent {\bf 1.}  In the same way as we can vary the size of the model in the $t$-direction (using the excess $\xi$)
we can vary the radius of the ball $\bB$. This can be done by reparametrizing the $x$-direction.
The reason we only chose to vary the $t$-direction is because these models are applied in \cite{O1}.

\noindent {\bf 2.} In the applications we may actually need warped metrics with warping functions that are multiples of hyperbolic
functions. All these functions are close to the exponential $e^t$ (for $t$ large), so instead of introducing one model for each hyperbolic function
we introduce only the exponential model.

Let $|.|_{C^k}$ denote the uniform $C^k$-norm of $\R^l$-valued functions on $\T_\xi=\bB\times I_\xi\sbs\R^{n+1}$. Sometimes we will write $|.|=|.|_{C^2}$.
Given a metric $g$ on $\T$, $|g|_{C^k}$ is computed considering $g$ as the $\R^{(n+1)^2}$-valued function $(x,t)\mapsto (g_{ij}(x,t))$ where, as usual,
$g_{ij}=g(e_i,e_j)$, and the $e_i$'s are the canonical vectors in $\R^{n+1}$. We will say that a metric $g$ on $\T$ is $\epsilon$-{\it close to hyperbolic}
if $|g-\sigma|_{C^2}<\epsilon$.\vspace{.1in}

A Riemannian manifold $(M,g)$ is $\epsilon$-{\it close to hyperbolic} if there is $\xi\geq 0$ such that for every $p\in M$ there is an {\it $\epsilon$-close to hyperbolic 
chart with center $p$}, that is, there is a chart
$\phi :\T_\xi\ra M$, $\phi(0,0)=p$,  such that $\phi^*g$ is $\epsilon$-close to hyperbolic. Note that all charts are defined on the same model space
$\T_\xi$. The number $\xi$ is called the {\it excess } of the charts
(which  is fixed).
More generally, a subset $S\sbs M$ is $\epsilon$-close to hyperbolic if every $p\in S$ is the center of an $\epsilon$-close to hyperbolic chart in $M$ with fixed
excess $\xi$. \vspace{.1in}

\noindent {\bf Remark.} 
Note that $\HH^n$ is $\epsilon$-close to hyperbolic for every $\epsilon >0$, but if a
hyperbolic manifold is not complete or have very small injectivity radius then it is not  $\epsilon$-close to hyperbolic.\vspace{.1in}

Let $I$ be an interval and consider $M\times I$ with variable metric $h=g_t+dt^2$. We say that a subset $S$ of $(M\times I ,h)$ is {\it radially} $\epsilon$-{\it close to hyperbolic}
if every $p\in S$ is the center of a {\it radially $\epsilon$-close to hyperbolic  chart of fixed excess $\xi$}, that is, by charts
$\phi :\T_\xi\ra M$, $\phi(0,0)=p$,
 such that $\phi^*g$ is $\epsilon$-close to hyperbolic and $\phi$ respects the ``structure" of $h$, i.e. $\phi$ satisfies the following two
conditions:
\begin{enumerate}
\item[{\bf (i)}] the map $\phi$ respects the $\B$-direction, i.e. $\phi \,(\bB\times\{ t\})\sbs
M\times\{ t'\}$
\item[{\bf (ii)}] the map $\phi$ preserves the $t$-direction, i.e. $\phi\, (\{ x\}\times I_\xi)\sbs
\{ x'\}\times I$ and $\phi |_{(\{ x\}\times I_\xi)}$ is an isometry for all $x\in\B$.
\end{enumerate}
Equivalently, the chart has the form $\phi(x,t)=(\phi_1(x),\, t+a)$, for some fixed number $a\in\R$ and chart $\phi_1$ on $M$. 
Note that, in this case, the pullback metric has the form $\phi^*h=\phi^*g_t+dt^2$, i.e. it is a variable metric on $\T$.
Note also that every radially $\epsilon$-close to hyperbolic chart $\phi$ can be extended to $\phi:\B\times (I-\{ a\})$ by the same formula $\phi (x,t)=(\phi_1(x),t+a)$,
but this extension may fail to be $\epsilon$-close to hyperbolic. (Here $I-\{a\}=\{t-a,\, t\in I\}$.)
Of course a  radially $\epsilon$-close to hyperbolic manifold is $\epsilon$-close to hyperbolic. Note that the definition of a radially $\epsilon$-close to hyperbolic metric depends on the product decomposition $M\times I$.\vspace{.1in}

\noindent {\bf Example.} Consider $M^{n+1}=\bS^n\times\R^+$ with warped metric $sinh^2t\,\sigma\0{\bS^n}+dt^2$. Hence $M$ is isometric to
a punctured hyperbolic space. As stated in the remark above $M$ is not $\epsilon$-close to hyperbolic,
but $S_L=\bS^n\times (L,\infty)\sbs M$ is $\epsilon$-close to hyperbolic, provided $L>2$. On the other hand $S_L$ is not radially $\epsilon$-close to hyperbolic 
when $\epsilon$ is small.
But the following is true: given $\epsilon>0$ there is $L>0$ such that $S_l$  is $\epsilon$-close to hyperbolic, for $l\geq L$. Actually, a more general 
statement is proven in 3.9 and 3.10.\vspace{.1in}

\noindent {\bf Remark 2.1.} 
For every $n$ there is a function $\epsilon'=\epsilon'(\epsilon,\xi)$ with the following  property:  if a Riemannian metric $g$ on a 
manifold $M^{n+1}$ is $\epsilon'$-close to hyperbolic, with charts of excess $\xi$,  then the sectional curvatures of $g$ all lie $\epsilon$-close to -1.
This choice is possible, and depends only on $n$ and $\xi$, because the curvature depends only of the derivatives up to order 2 of $\phi^*g$ on $\T_\xi$,
where $\phi$ is an $\epsilon$-close to hyperbolic chart with excess $\xi$.\vspace{.3in}

\noindent {\bf \large 3. Slow families of metrics and the proof of the Main Theorem.}

Consider the family of metrics $g\0{t}$, $t\in I\sbs \R$, on $M$, $M$ closed. We say that 
$\{ g\0{t}\}$ is {\it  $\epsilon$-slow} if $g\0{t}$ and its first derivatives 
change $\epsilon$-slowly with
$t$. 
 That is:
 \begin{enumerate}
\item[(i)] for every $t_0\in I$, \,$k=1,\,2$ and \,$u\in TM$ we have
 $\Big|\frac{d^k}{d\,t^k}\,g\0{t}(u,u)|\0{t_0}\,\Big|\,\leq\,\epsilon\,\, g\0{t_0}(u,u)$
\item[(ii)] for every $t_0\in I$, $v\in TM$ and   $u$ vector field on $M$ we have 
\end{enumerate}
$$\Big|\,\frac{d}{d\,t}\,v\,g\0{t}(u,u)|\0{t_0}\,\Big|\,\leq\,\epsilon\,\Big( \,g\0{t_0}(u,u)\,
g\0{t_0}^{{\mbox{\tiny 1/2}}}(v,v)\,+\,
g\0{t_0}^{{\mbox{\tiny 1/2}}}(u,u)\,\,g\0{t_0}^{{\mbox{\tiny 1/2}}}(\n _vu,\n_vu)\,\Big)$$
For instance, if the family $\{ g\0{t}\}$ is constant, then it is $\epsilon$-slow, for every $\epsilon>0$. We will  need the following lemma later.\vspace{.1in}

\noindent {\bf Lemma 3.1.} {\it Consider the family of metrics $\{ g\0{t}\}_{t\in I}$, $I\sbs \R$ an interval, and assume
$\{ g\0{t}\}$ is $\epsilon$-slow. 

\begin{enumerate}
\item[(1)] The family of metrics $\{ g\0{t+b}\}_{t\in I-b}$ is $\epsilon$-slow, for any $b\in\R$. 

\item[(2)] Let $\varphi: J\ra I$ be a diffeomorphism with $|\varphi'(t)|,\,|\varphi''(t)|<a$,
$t\in J$. Then
the family $\{ g_{\varphi(s)}\}_{s\in J}$ is $\big(\,(a+a^2)\epsilon)\,\big)$-slow.

\item[(3)] Let $\phi:M\ra M$ be a diffeomorphism. Then
 $\{\phi^* g\0{t}\}_{t\in I}$ is $\epsilon$-slow.\end{enumerate}}

\noindent{\bf Proof.} Statement (1) is direct. Statment (2)
follows from the chain rule, and (3) follows from the fact that
the definition of $\epsilon$-slow metrics is intrinsic. This
proves the lemma.\vspace{.1in}

Let $c>1$.
A metric $g$ on a compact manifold $M$ is {\it c-bounded} if $|g|< c$ and $|\, det \,g\,|_{C^0} > 1/c$.
A set of metrics $\{ g\0{\lambda}\}$ on the compact manifold $M$ is {\it c-bounded} if every $g\0{\lambda}$
is $c$-bounded.  \vspace{.1in}

\noindent {\bf Remarks.}

\noindent {\bf 1.} Here the uniform $C^2$-norm $| .|$ is taken with respect to a fixed finite atlas $\cA$. Hence the definition of
a $c$-bounded family depends on the choice of the atlas $\cA$. 

\noindent {\bf 2.} We will assume that the finite atlas $\cA$ is
``nice", that is: (1) it has``extendable" charts, i.e.
charts that can be extended to the (compact) closure of their domains, and (2) we assume the $\cA$ has the following
property: for every $p\in M$ there is a a chart $\phi:U\ra M$ in $\cA$
with $d\0{\R^n}(\phi^{-1}(p), \R^n-U)\geq 1$.
If an atlas does not have this property then each chart
can be precomposed with a dilation to obtain an new atlas with
this desired property. 

\noindent {\bf 3.} For a metric $g$ on $M$ we will use the same symbol $g$ to denote a matrix representation in a chart of the
fixed finite atlas.

\noindent {\bf 4.} We are taking $c >1$ just to simplify  some calculations.

\noindent {\bf 5.} If $\{ g\0{t}\}_{t\in I}$ is a (smooth) family and $I\sbs \R$ is compact then, by compacity and continuity, the family
 $\{ g\0{t}\}$ is bounded. (Recall we are taking $M$ compact.)

\noindent {\bf 6.} If $\{ g\0{t}\}_{t\in I}$ is $c$-bounded then clearly
$\{ g\0{t(s)}\}_{s\in J}$ is also $c$-bounded, for any reindexation (or reparametrization)
$t=t(s)$.\vspace{.1in}

For an excess $\xi$ and $I\sbs\R$ we write $I(\xi)=\big\{\,  t\in I,\,  \big(\,t-(1+\xi), t+(1+\xi)\,\big)\sbs I\,\big\}$; equivalently $I(\xi)$ is the maximal set such that 
$I(\xi)+I_\xi\sbs I$. We write 
$T=inf \,I(\xi)=1+\xi+inf\, I$. Hence $t\geq T$, for all $t\in I(\xi)$.  
We will assume $inf\, I$ to be positive, therefore we will always have $T> 1$.
The following Theorem is slightly different from the Main Theorem: 
it considers the excess of the charts and different warping function. This is the result will prove.\\

\noindent {\bf Theorem 3.2.} {\it Let $M^n$ be a closed smooth manifold
and $\xi\geq 0$.
If the family of metrics $\{g\0{t}\}_{t\in I}$ on $M$ is $\epsilon$-slow and $c$-bounded then
the metric $h=e^{2t}g\0{t}+dt^2$ is radially $\eta$-close to hyperbolic on $M\times I(\xi)$
with charts of excess $\xi$, provided $$C\,(e^{-T}+\epsilon)\leq\eta$$ 
\noindent where $C=C(c,n, \xi)$.}\\

\noindent {\bf Remarks.} 

\noindent {\bf 1.} The constant $C$ depends solely on the dimension $n$
of $M$, the constant $c$ and the desired excess $\xi$. An explicit formula for $C$ is given
at the end of the proof of Theorem 3.2. Note that if $T$ is large and $\epsilon$ small
then we can take $\eta$ small.

\noindent {\bf 2.} Recall that the definition of a $c$-bounded family
depends of the choice of the atlas $\cA$, and we are assuming
that $\cA$ is ``nice" (see Remarks 1, 2 after 3.1). But if $\cA$ is
not nice the constant $C$ in Theorem 3.2 would depend on a certain Lebesgue
type number of the atlas $\cA$.

\noindent {\bf 3.} The construction of the $\eta$-close to hyperbolic
charts around a given point is explicit; see formula $(\star)$ in the
proof of Theorem 3.2 below.

\noindent {\bf 4.} Recall we are assuming $T> 1$. This is just for convenience;
things could be arranged to include the case $T\leq1$.
\vspace{.1in}

The following Corollary implies the Main Theorem by taking the excess $\xi=0$.\\

\noindent {\bf Corollary 3.3.} {\it Under the same conditions as in 3.2 we have that 
the metric $h=sinh^2t\,g\0{t}+dt^2$ is radially $\eta$-close to hyperbolic on $M\times I(\xi)$,
with charts of excess $\xi$, 
provided $$\,C_1\,(e^{-T}+\epsilon)\leq\eta$$ 
\noindent where $C_1=C_1(c,n,\xi)$.}\\

\noindent {\bf Remarks.}

\noindent {\bf 1.} We can take $C_1=30\, C(6^{n}c,n,\xi)$, where
$C$ is as in Theorem 3.2.

\noindent {\bf 2.} It can easily be checked from the proofs of 3.2 and 3.3 that the  $\eta$-close to hyperbolic
charts in Corollary 3.3 are the same as the ones
constructed in the proof of Theorem 3.2; see formula $(\star)$ in the
proof of 3.2, and also see Remark 2 after 3.2 above.

\noindent {\bf 3.} Recall we are assuming $T>1$ (see Remark 4 above).
\vspace{.1in}

We will use the following lemmas.\vspace{.1in}

\noindent {\bf Lemma 3.4.} {\it Let $G$ be a real positive definite symmetric $n\times n$ matrix with $|G_{ij}|<c$ and $|det\,G|>1/c$. Then we can write $G= F^TF$, with $|F_{ij}|\,<n\sqrt{c}$ and $ |(F^{-1})_{ij}|<\, n\sqrt{(n-1)!\,c^{n}}$.}\vspace{.1in}

\noindent {\bf Proof.} Since $G$ is symmetric we can write $G=O^T D O$, where the columns of $O^{T}$ form an orthonormal
basis (in $\R^n$) of eigenvectors of $G$ and $D$ is a diagonal matrix whose diagonal entries are the corresponding eigenvalues $\mu$ of $G$. 
Since $D=OG O^T$, $O$ is orthonormal, $|G_{ij}|<c$ and $|det \,G| > 1/c$,   we have that 
$ |\mu |\, <\, n c$. On the other hand, since $|det\, G|>1/c$,  we have that $|(G^{-1})_{ij}|<(n-1)! c^n$, and a similar argument as before applied to $D^{-1}$
shows $|1/\mu|\,<\,n!\,c^{n}\,$. Hence 
$\frac{1}{n!\,c^{n}} \, <|\mu |\, <\, n c$.
Take $F=\sqrt{D}O$. This proves the lemma.\vspace{.1in}

\noindent {\bf Remark.} Let $u\in\R^n$, and let $|u|$ be the euclidean norm.
The following estimate follows from the proof of 3.4:
$$
\frac{1}{n!\, \,c^{n}}\,\,|u|^2\,<\,| u^TGu|\,<n\,c\, |u|^2
$$

\noindent {\bf Lemma 3.5.} {\it  Let $\{ g\0{t}\}$ be an $\epsilon$-slow $c$-bounded family of metrics on the unit ball $\bB^n\sbs\R^n$.
Then }
\begin{enumerate}
\item[]  $|\frac{\p^k}{\p t^k}\big( g\0{t}\big)\0{ij}|\0{C^0}\leq 3\,\epsilon \, c $,\,\, $k=1,\,2$.
\item[]  $|\frac{\p^2}{\p t\,\p x\0{l}}\big( g\0{t}\big)\0{ij}|\0{C^0}\leq \epsilon \, c\0{3} $
\end{enumerate}
\noindent {\it where $c\0{3}=c\0{3}(c,n)=3\,c^{3/2}+\frac{9}{2}n!
c^{n+2}$.}\vspace{.1in}

\noindent {\bf Proof.} We prove the first inequality for $k=1$. The proof for $k=2$
is similar. To simplify our notation write $g'\0{ij}=\frac{\p}{\p t}(g\0{t})\0{ij}$.
Also sometimes we will omit the variable $t$.
From the definition of slow family metrics we have that 
$|\frac{\p}{\p t} g\0{t}(u,u)(t_0)|\leq \epsilon \,  g\0{t_0}(u,u)$. 
In what follows of this proof the bars $|.|$ denote the $C^0$ norm. Taking $u=\p_i$ 
and using the fact that the family is $c$-bounded we get that $|g'\0{ii}|\leq\epsilon\, c$.
And taking $u=\p_i+\p_j$ we get $|g'\0{ii}+2g'\0{ij}+g'\0{jj}|\leq
\epsilon |g\0{ii}+2g\0{ij}+g\0{jj}|\leq 4\epsilon\,c$. These two inequalities imply
$$
|2g'\0{ij}|\leq |g'\0{ii}+2g'\0{ij}+g'\0{jj}|+|g'\0{ii}|+|g'\0{jj}|\leq\,6\,\epsilon\,c
$$
\noindent This proves the first inequality.  We prove the second inequality. Since
the family of metrics is $c$-bounded we have that $|g^{ij}|<(n-1)!\,c^{n}$, where $(g^{ij})$ is the inverse of $g=(g\0{ij})$. Hence we obtain the following estimate for
the Christoffel symbols

\begin{equation*} \Big|\Gamma^k\0{ij}\Big|\,\,<\,\,\frac{3}{2} \,(n-1)!\, c^{n+1} 
\tag{a}
\end{equation*}

\noindent In what follows of this proof we use the summation notation. It
follows from (a) above that

{\small \begin{equation*}
g\Big( \n_{\p_l}\p_i\,,\, \n_{\p_l}\p_j \Big)\,\,=\,\, g\Big( \Gamma\0{li}^k\p_k,
\Gamma\0{lj}^m\p_m \Big)\,\,=\,\, g\0{km}\Gamma\0{li}^k\Gamma\0{lj}^m
\,\,<\,\, \frac{9}{4}\,\big(n!\big)^2\,c^{2n+3}\,=\, c\0{1}^2
\tag{b}
\end{equation*}}

\noindent Write $g'\0{ij,l}=\frac{\p^2}{\p t\,\p x\0{l}} g\0{ij}$.
Now take $v=\p_l$ and $u=\p_i$ in (ii) of the definition of $\epsilon$-slow
metrics, and use (b) to obtain
\begin{equation*}
\Big| g'\0{ii,l}\Big|\,\leq\,\epsilon\, \Big(\, c\, c^{1/2}\, +\,
c^{1/2}\, c\0{1}\,  \Big)\,=\, \epsilon\, c\0{2}
\tag{c}
\end{equation*}
\noindent where $c\0{2}=c\0{2}(c,n)= c^{3/2}+c^{1/2}\,c\0{1}=c^{3/2}+
\frac{3}{2}\,n!\,c^{n+2}$.
Take $v=\p_l$,  $u=\p_i+\p_j$ now in (ii) in the definition of $\epsilon$-slow
metrics and a calculation using (b), (c) show
$$\begin{array}{lll}
\Big| 2g'\0{ij,l} \Big|&\leq&
\Big| g'\0{ii,l}+2g'\0{ij,l}+g'\0{jj,l} \Big|\,+\,\Big| g'\0{ii,l} \Big|\,+\,
\Big| g'\0{jj,l} \Big|\\  \\&\leq& \epsilon\,\Big(\,(4\,c)\,c^{1/2}\,+\,2\, c^{1/2}\,(2\,c\0{1}) \Big)
\,+\,\epsilon \,c\0{2}\,+\,\epsilon \,c\0{2}
\,\,=\,\, 2\,\epsilon\, c\0{3}
\end{array}
$$

\noindent This proves the lemma.\vspace{.1in}

We shall prove a sort of a converse to this lemma in 3.6.\vspace{.1in}

\noindent {\bf Proof of Theorem 3.2.} Assume the family of metrics
$\{g\0{t}\}_{t\in I}$ is $c$-bounded and $\epsilon$-slow. Recall that we are denoting by $\bB_a(x)\sbs\R^n$ the ball of
radius $a$ centered at $x$.\vspace{.1in}

First we reduce the problem to $\R^n$ using the fixed finite atlas $\cA$: 
for each $p\in M$ choose a chart $\psi_p=(\varphi,\bB_1(x))$, with $\varphi(x)=p$, such that $\psi_p$ is the restriction of one of the charts in $\cA$ (see Remarks 1 and 2 after 3.1).
And after a  translation we can assume  $x=0$, that is,
$\varphi:\bB=\bB_1(0)\ra M$, $\varphi(0)=p$.\vspace{.1in}

Fix $\xi\geq 0$, $p\in M$ and $t_0\in I(\xi)$. Thus $(t_0-(1+\xi),t_0+(1+\xi))\sbs I$. Write $g\0{t_0}(p)=F^TF$ with $F$ as in Lemma 3.4.
(Here $g\0{t_0}$ denotes also the matrix representation of $g\0{t_0}$ in the chart $\psi_p$.) Let $A=F^{-1}$.
By Lemma 3.4 we have  $|A_{ij}|< c\0{4}$, where $c\0{4}=\sqrt{n\,n!\, c^{n}}$. 
Define the chart $\phi :\T=\bB\times I_\xi\ra \bB\times I\sbs M\times I$ by
\begin{equation*} \phi (x,t)=\Big(\,\varphi\big(\,e^{\lambda-t_0}A x\,\big)\,, t+t_0\,\Big)\tag{$\star$}
\end{equation*}

\noindent where $\lambda=min\{0, t_0-ln(nc\0{4})\}$.
The proof of Theorem 3.2 consists of showing that the chart
$\phi$ is $\eta$-close to hyperbolic, with $\eta\geq C(e^{-t_0}+\epsilon)$,
for certain $C=C(c,n,\xi)$. (Recall $t_0\geq T$.) To do this we need to estimate $|f-\sigma|$ (this is the $C^2$ norm, and recall that $\sigma$ is the model metric, see 
Section 2). But before, in the next claim,  we show that $\phi$ is defined
on the unit ball $\B$:\vspace{.1in}

\noindent {\it Claim.} {The map $x\mapsto e^{\lambda-t_0}A x$
sends $\bB$ into itself.}\vspace{.1in}

\noindent This claim was the reason for introducing $\lambda$, which is a correction term for $t$ small. \vspace{.1in}

\noindent {\it Proof of Claim.}
We have $|A|<c\0{4}$ (here $|.|$ is the uniform
norm). Therefore $|Ax|<nc\0{4}$ for $x\in\bB$ (here
$|.|$ is the Euclidean norm). Therefore $|e^{\lambda-t_0}Ax|<
e^{\lambda-t_0}nc\0{4}$ for $x\in\bB$. We have two possibilities: (1) $t_0\geq
ln(nc\0{4})$ or (2) $t_0\leq ln(nc\0{4})$. In case (1) we have
$\lambda=0$ and $e^{\lambda-t_0}=e^{-t_0}<1/nc\0{4}$, hence
$|e^{-t_0}Ax|<1$, for $x\in \bB$. In case (2) we have $\lambda=
t_0-ln(nc\0{4})$, thus $\lambda-t_0=-ln(nc\0{4})$ and also follows that
$|e^{\lambda-t_0}Ax|<1$, for $x\in \bB$. This proves the claim.
\vspace{.1in}

Let $h=e^{2t}g_t +dt^2$ as in the statement of 3.2.
Write $f=\phi^*h$.  As mentioned before we want to estimate $|f-\sigma|$.\vspace{.1in}

Since $f=\phi^*h$ we have $f=e^{2t}f_t+dt^2$, where $f_t=e^{2t_0}\,\varphi^*g\0{t+t_0}$,
$t\in I_\xi$. Then 
\begin{equation*} f_t(x)=e^{2\lambda}A^T\, g\0{t+t_0}(e^{\lambda-t_0}Ax)\, A\tag{1}\hspace{.4in}{\mbox{and}}\hspace{.3in}(f_t)\0{ij}=e^{2\lambda}\sum\0{k,l} A\0{li}g\0{lk}A\0{kj}
\end{equation*}

\noindent By hypothesis we have $|g\0{ij}|<c$ (recall this is $C^2$ norm).
Evaluating $f$ at $(0,0)$ we get that
\begin{equation*} f(0,0)=\phi^*h(0,t_0)=e^{2\lambda}\sigma_{\R^n}
+dt^2\tag{2}
\end{equation*}

\noindent (Note that $\sigma(0,0)=\sigma\0{\R^n}+dt^2$). Differentiating $f=e^{2t}f_t+dt^2$ and using equation (1)  and the facts that
$|A_{ij}|< c\0{4} $ and $e^{\lambda}\leq 1$ we get the following estimates:
\begin{equation*}  |\p_J f|\0{C^0}=e^{2t}|\p_J f_t|\0{C^0}\leq e^{2t}\, n^2\, c^2\0{4}\, \Big[ n\, c\,
e^{\lambda-t_0}c\0{4}\Big]^{|J|} \leq e^{2t}\, n^2\, c^2\0{4}\, \Big[ n\, c\, e^{-t_0}c\0{4}\Big]^{|J|}\tag{3}\end{equation*}
\noindent where $J$ is a multi-index of order $|J|=1,\,2$ in the $\bB$-direction, i.e. no $t$-derivatives are considered.
It follows from Lemma 3.5 (applied to the family $\{g_t\}$)
and (1) that (recall $e^\lambda\leq 1$)

\begin{equation*} |\frac{\p^k}{\p t^k}\big( f_t\big)_{ij}|\0{C^0}\leq 3\,n^2\,c^2\0{4}\,c\,  \epsilon\,=\,c\0{5}\tag{4}\end{equation*}

\begin{equation*} |\frac{\p^{2}}{ \p x\0{l}\p t}\big( f_t\big)_{ij}|\0{C^0}\leq \epsilon \,\bigg[ n^3\,
c\0{3}\,c\0{4}^3\, e^{-t_0}\bigg]\,=\,c\0{6} 
\tag{5}\end{equation*}

\noindent Recall that the matrix representation
 of $\sigma\0{\R^n}$ is the identity matrix $1\!\!1$. Thus
 $|f_t(0)-\sigma\0{\R^n}|=|f_t(0)-1\!\!1|$, $t\in I_\xi$.
 And from (2) and (4), for $t\in I_\xi$  we get
\begin{equation*} |f_t(0)-1\!\!1|\leq|f_0(0)-1\!\!1|+\int_0^{1+\xi}|\frac{\p}{\p t}f_t(0)| dt\,\leq\, (1-e^{2\lambda})\,+
\, (1+\xi)\,c\0{5} 
\end{equation*}

\noindent This together with $|1-e^{2\lambda}|=[(1-e^{2\lambda})e^{t_0}]e^{-t_0}\leq nc\0{4} e^{-t_0}$ (recall $\lambda=0$ when
$e^{t_0}\geq nc\0{4}$) imply

\begin{equation*} |f_t(0)-1\!\!1|\leq\, nc\0{4} e^{-t_0}\,+
\, (1+\xi)\,c\0{5} \,=\,c\0{7}
\tag{6}\end{equation*}

\noindent Write $c\0{8}=n^3\,c\0{4}^3\,c\,e^{-t_0}$. Then from (3) we have
$|\frac{\p}{\p x\0{l}}f_t|\0{C^0}\leq c\0{8}$.
This together with (6) and the Mean Value Theorem imply
$$
|f(x,t)-\sigma(x,t)|\,\leq\,e^{2t}\bigg( |f_t(x)-f_t(0)|+|f_t(0)-1\!\!1 | \bigg)\,\leq\,
e^{2(1+\xi)}\Big( n^{1/2}\,c\0{8}\,+c\0{7} \Big)\,=\, c\0{9}
$$
\noindent Hence
\begin{equation*} |f-\sigma|_{C^0}\,\leq \, c\0{9}
\tag{7}\end{equation*}

\noindent Now, since the derivatives in the $\bB$-direction of $\sigma$ all vanish,  equation (3)
holds replacing $f$ by $f-\sigma$ and we get
\begin{equation*} \Big| \,\p_J(f-\sigma)\,\Big|\0{C^0}\,\leq\,e^{2(1+\xi)}\, n^4\, c^4\0{4}\, c^2\,
e^{-t_0}\,=\, c\0{10}
\tag{8}\end{equation*}
\noindent where we are assuming $t_0>0$ so that $e^{-2t_0}<e^{-t_0}$.
Now, note also that $\pt (f_t-1\!\!1 )=\pt f_t$. This together with
the definitions of $f$ and $\sigma$ imply
$\pt (f-\sigma)=2(f-\sigma)+ e^{2t}\pt f_t$, which together with (4) and (7)
imply 
\begin{equation*} |\pt (f-\sigma)|\0{C^0}\,\leq\, 2\,c\0{9}\,+\, e^{2(1+\xi)}\, c\0{5}\,=\, c\0{11}
\tag{9}
\end{equation*}

\noindent Analogously, differentiating with respect to $t$ again and using (4) we get
\begin{equation*} |\frac{\p^2}{\p t^2} (f-\sigma)|\0{C^0}\,\leq\, 4\,c\0{9}\,+\, 4\, e^{2(1+\xi)}\,c\0{5}\,+\, e^{2(1+\xi)}\,c\0{5}\,=\, 4\,c\0{9}\,+\, 5\,e^{2(1+\xi)}\,c\0{5}\,=\, c\0{12}
\tag{10}
\end{equation*}
\noindent And from (5) and (8) we finally get

\begin{equation*} |\frac{\p^2}{\p x\0{l} \p t} (f-\sigma)|\0{C^0}\, \leq\, 2\,c\0{10}\,+\, e^{2(1+\xi)}\, c\0{6} \,=\,c\0{13}
\tag{11}
\end{equation*}

\noindent Hence $|f-\sigma|<c\0{14}$, where $c\0{14}=$max$\{
c\0{9},\,c\0{10},\,c\0{11},\,c\0{12},\,c\0{13} \}\leq C\big(\,e^{-t_0}\,+\,\epsilon\big)$
(recall the bars $|.|$ denote the $C^2$ norm). But each term in each of the constants
$c\0{i}$, $i=9,...,13$ contains an $\epsilon$ or $e^{-t_0}$ (or both).
And it can be verified that we can take
(assuming $n\geq 2$ and $t_0>0$)

$$C\,=\,C(c,n,\xi)\,=\,(27+4\xi)\,e^{2(1+\xi)}\, n^4\, c\0{4}^4\, c\0{3}\,c^2$$

\noindent (This is by no means an optimal choice). This completes the proof of Theorem 3.2.\vspace{.1in}

\noindent {\bf Proof Corollary 3.3.}
Let $k(t)=\frac{sinh^2(t)}{e^{2t}}=\frac{(1-e^{-2t})^2}{4}$.
A quick calculation shows that for $t\geq 1$ we have 
(note that we indeed have $t\geq1$ because $\xi\geq0$)

\begin{equation*}\begin{array}{l}
\frac{1}{6}\,<k\,<\frac{1}{4}\\\\
 0<\,k'\,<\,e^{-2t}\,<\,1/2\\\\
0<\,|k''|\,<\,2\,e^{-2t}\,<\,1
\end{array}
\tag{12}
\end{equation*}

We can write $h=sinh^2(t)\,g\0{t}+dt^2=e^{2t}\,k\,g\0{t}+dt^2$. Then the Corollary
follows from Theorem 3.2 and the following claim.\vspace{.1in}

\noindent {\bf Claim.} {\it The family of metrics $\{k\,g_t\}$ is
$30\,(e^{-2T}+\epsilon)$-slow and $(6^{n}\,c)$-bounded.
(Recall, by hypothesis, we have $t\geq T$).}

\noindent {\bf Proof of claim.} Using (12) we have 
$$|k\,g\0{t}|\0{C^0}=k\,|g\0{t}|\0{C^0}<\frac{1}{4}\,c\,<\, c$$
\noindent Also, $|\, det\,(k\,g\0{t})\,|\0{C^0}\,=\,k^n\,|\,det\,g\0{t}\,|\0{C^0}\,>\,\frac{1}{6^n}\,|\,det\,g\0{t}\,|\0{C^0}\,>\,\frac{1}{6^n\,c}$.
Hence $\{g\0{t}\}$ is $(6^{n}\,c)$-bounded. We prove now the first statement of the definition of $\epsilon$-slow metrics.
Using (12) and the fact that $\{g\0{t}\}$ is $\epsilon$-slow we get
\begin{equation*} \begin{array}{lllll}\Big( (k\, g)(u,u)\Big)'&=&k'\,g(u,u)\,+\, k\,g'(u,u)
\,\,\,\,\leq\,\,\,\,
\frac{1}{k}\,\big(k'\,+\, k\,\epsilon\big)\,(k\,g)(u,u)\\\\
&\leq&6\,\big(\,e^{-2T}+\epsilon\,\big)\,(k\,g)(u,u)\,\,\,\,\leq\,\,\,\,
30\,\big(e^{-2T} +\epsilon\big)\,(k\,g)(u,u)\end{array}\tag{13}\end{equation*}
\noindent and
\begin{equation*} {\small
\begin{array}{lllll}\Big( k\, g(u,u)\Big)''&=&k''\,g(u,u)\,+\, 2\,k'\,g'(u,u)\,+\,k\,g''(u,u)
\,\,\,\,\leq\,\,\,\,
\frac{1}{k}\,\big(k''\,+\,2\, k'\,\epsilon\,+\, k\,\epsilon\big)\,(k\,g)(u,u) \\ \\
&\leq &6\,\big(\,2\,e^{-2T}\,+\,\epsilon\,+\, \epsilon\,\big)\,(k\,g)(u,u)\,\,\,\,\leq \,\,\,\,30\,\big(\,e^{-2T} +\epsilon\,\big)\,(k\,g)(u,u)
\end{array}}
\tag{14}\end{equation*}
\noindent Also
\begin{equation*}{\small \begin{array}{lll}
\Big(v\,k\,g(u,u)\Big)'&=&k'\Big(v\,g(u,u)\Big)\,+\, k\,\Big(v\,g(u,u)\Big)'\\\\
&\leq &k'\Big(v\,g(u,u)\Big)\,+\,\epsilon\,k\,\Big(\,g(u,u)
\,g^{^{\mbox{\tiny $1/2$}}}(v,v)\,+\,
\,g^{^{\mbox{\tiny $1/2$}}}(u,u)\, g^{^{\mbox{\tiny $1/2$}}}
(\n_vu,\n_vu)   \Big)
\end{array}}
\tag{15}
\end{equation*}
\noindent But
\begin{equation*}
\Big|\,v\,g(u,u)\,\Big|\,=\,2\,\Big|\,g(u,\n_vu)\,\Big|\,\leq\,
2\,g^{^{\mbox{\tiny $1/2$}}}(u,u)\,\,g^{^{\mbox{\tiny $1/2$}}}(\n_vu,\n_vu)
\tag{16}
\end{equation*}
\noindent And from (15),(16) and (12) we get
\begin{equation*}{\small\begin{array}{lllllll}
\Big(v\,k\,g(u,u)\Big)'&\leq &\Big(\,2\,k'\,+\,k\,\epsilon\Big)
\Big(\,g(u,u)
\,g^{^{\mbox{\tiny $1/2$}}}(v,v)\,+\,
\,g^{^{\mbox{\tiny $1/2$}}}(u,u)\, g^{^{\mbox{\tiny $1/2$}}}
(\n_vu,\n_vu)   \Big)&&&\\\\&\leq &
\Big(\,2\,k'\,+\,k\,\epsilon\Big)
\Big(\,g(u,u)
\,g^{^{\mbox{\tiny $1/2$}}}(v,v)\,+\, k^{-1/2}
\,g^{^{\mbox{\tiny $1/2$}}}(u,u)\, g^{^{\mbox{\tiny $1/2$}}}
(\n_vu,\n_vu)   \Big)&&&\\\\
&\leq &k^{-3/2}\,
\Big(\,2\,k'\,+\,k\,\epsilon\Big)
\Big(\,k^{3/2}\,g(u,u)
\,g^{^{\mbox{\tiny $1/2$}}}(v,v)\,+\, k
\,g^{^{\mbox{\tiny $1/2$}}}(u,u)\, g^{^{\mbox{\tiny $1/2$}}}
(\n_vu,\n_vu)   \Big)&&&\\\\
&=&15\,
\Big(2\,k'+k\,\epsilon\Big)
\Big((k\,g)(u,u)
\,(k\,g)^{^{\mbox{\tiny $1/2$}}}(v,v)+ 
\,(k\,g)^{^{\mbox{\tiny $1/2$}}}(u,u)\, (k\,g)^{^{\mbox{\tiny $1/2$}}}
(\n_vu,\n_vu)   \Big)&&&
\\\\&\leq &30\,\big(\,e^{-2T} +\epsilon\,\big)\,
\Big((k\,g)(u,u)
\,(k\,g)^{^{\mbox{\tiny $1/2$}}}(v,v)+ 
\,(k\,g)^{^{\mbox{\tiny $1/2$}}}(u,u)\, (k\,g)^{^{\mbox{\tiny $1/2$}}}
(\n_vu,\n_vu)   \Big)&(17)&&
\end{array}}
\end{equation*}

\noindent And the claim follows from (13), (14) and (17). This proves the claim.\vspace{.1in}

The Corollary now follows from Theorem 3.2, the claim and the fact that
(recall $T>1$) $$
e^{-T}\,+\,30\,(\,e^{-2T}\,+\, \epsilon\,)\, \leq \,30\,(\,e^{-T}\,+\,\epsilon\,)
$$
\noindent This proves Corollary 3.3.\vspace{.1in}

Here is a simple particular case of Theorem 3.2.  Assume that the family $\{ g\0{t}\}$ is constant, i.e. $g\0{t}=g\0{0}$, for all $t\in I$, 
then the variable metric
$g=g\0{t}+dt^2=g\0{0}+dt^2$ is just a product metric, and the  metric $h=e^{2t}g\0{0}+dt^2$ is a warped metric.
In this case Theorem 3.2 says that, given $\epsilon >0$, the warped metric $h$ is $\epsilon$-close to hyperbolic, provided $t$ is large
enough (how large depending on $\epsilon$, the dimension $n$ and the metric $g\0{0}$).\vspace{.1in}

Here is a sort of a converse to Lemma 3.5.\vspace{.1in}

\noindent {\bf Lemma 3.6.} {\it Suppose that the family $\{g\0{t}\}$
(on some open $U\sbs\R^n$)  is $c$-bounded and that}
\begin{enumerate}
\item[]  $|\frac{\p^k}{\p t^k}\big( g\0{t}\big)\0{ij}|\0{C^0}\leq \epsilon\0{1} $,\,\, $k=1,\,2$.
\item[]  $|\frac{\p^2}{\p t\,\p x\0{l}}\big( g\0{t}\big)\0{ij}|\0{C^0}\leq \epsilon\0{2} $
\end{enumerate}
\noindent {\it Then $\{ g\0{t}\}$ is $\epsilon\0{3}$-slow, where
$\epsilon\0{3}=\Big(\epsilon\0{2}\,n\,+\,2\,\epsilon\0{1}\, c\0{14}\,    \Big)\,n^2
\Big(n!\,c^{n}\Big)^{3/2}$, and $c\0{14}=\frac{3}{2}n^{3/2}n!c^{n+2}$.}

\noindent {\bf Remark.} Note that $\epsilon\0{3}=a\,\epsilon\0{1}\,+\,b\,\epsilon\0{2}$
where $a=a(n,c)$ and $b=b(n,c)$ are constants that depend solely on $c$ and $n$.

\noindent {\bf Proof.} 
First we prove (i) of the definition of $\epsilon$-slow metrics. We do this just for $k=1$
because the proof for $k=2$ is similar. We use the summation notation and denote
the derivative with respect to $t$ by a prime. Also the euclidean norm on $\R^n$
will be denoted by bars $|.|$. 

Write $v=v^i\p_i\in T_x\R^n$. Then using the remark after the proof of 3.4 we get 
$$g(v,v)'\,=\,(g\0{ij}v^iv^j)'\,=\,g\0{ij}'v^iv^j\,\leq\,\epsilon\0{1}n^2\,|v|^2\,<\,
\epsilon\0{1}\,n^2\,n!\, c^{n}\, g(v,v)$$
\noindent This proves (i). To prove (ii) first we show the following estimate: (here $v=v^k\p_k
\in T_x\R^n$ and $u=u^l\p_l$ is a vector field on $U$) \begin{equation*}
\big| v(u^k) \big|\,\leq\,\big| v(u^k)\p_k \big|\,\leq\, \big|  \n_vu \big|\,+\,c\0{14} \big|v \big|\,
\big| u_x \big|
\tag{d}
\end{equation*}

\noindent where $c\0{14}=\frac{3}{2}n^{3/2}n!c^{n+2}$ and $u_x$ is the value of $u$ at $x$. To prove this note that 
\begin{equation}\n_vu=\n_vu^k\p_k=v(u^k)\p_k+u^k_xv^l\n_{\p_l}\p_k=v(u^k)\p_k+
u^k_xv^l\Gamma^s_{kl}|_x\p_s\tag{e}\end{equation}
\noindent Equation (d) now follows from (e) and (a) in the proof of 3.5.

Now, to prove (ii) we use (d)  to compute:
\begin{center}$\begin{array}{lllll}
\Big| \frac{d}{dt}v\,g (u,u) \Big|&=&
\bigg[ v\,\Big( g\0{jk}u^ju^k  \Big)   \bigg]'& = &v^l\Big(\p\0{l}\,g\0{jk}\Big)'u^j_xu^k_x\,+
\,2\,g'\0{jk}\,u^j_x\,v(u^k)\\ &&&<&\epsilon\0{2}\, n^3\, |v|\,|u_x|^2\,+\,2\,
\epsilon\0{1}n^2\,|u_x|\,\Big(|\n_vu|+c\0{14}|u_x|\,|v|\Big)\\
&&&=&\Big(\epsilon\0{2}\,n+2\,\epsilon\0{1}\,c\0{14}\Big)\,n^2\,|v|\,|u_x|^2\,+\,
2\,\epsilon\0{1}\,n^2\,|u_x|\,|\n_vu|\\
\end{array}
$\end{center}
\noindent This together with the remark after 3.4 gives:
\begin{center}$\Big| \frac{d}{dt}v\,g(u,u) \Big|\,\,<\,\,
 A\,\,g\0{t_0}(u,u)\,
g\0{t_0}^{{\mbox{\tiny 1/2}}}(v,v)\,\,+\,\, B\,\,
g\0{t_0}^{{\mbox{\tiny 1/2}}}(u,u)\,\,g\0{t_0}^{{\mbox{\tiny 1/2}}}(\n _vu,\n_vu)
$\end{center}

\noindent where $A=\Big(\epsilon\0{2}\,n\,+\,2\,\epsilon\0{1}\, c\0{14}\,    \Big)\,n^2
\Big(n!\,c^{n}\Big)^{3/2}$ and $B=\Big(2\,\epsilon\0{1}\, n^2  \Big)\,
\Big( n!\, c^{n} \Big)$. Since $A>B$ the lemma follows.
This proves the lemma.\vspace{.1in}

Recall we are assuming $M$ closed. If $I$ is compact then we can find
$\epsilon\0{1}$ and $\epsilon\0{2}$ as in  the statement of Lemma 3.6.
Hence we obtain the following Corollary.\vspace{.1in}

Here is a sort of a local converse of Theorem 3.2.\vspace{.1in}

\noindent {\bf Corollary 3.7.} {\it Let $g=e^{2t}g\0{t}+dt^2$ be a variable metric on $\T_\xi$. Assume $g$ is $\epsilon$-close to hyperbolic. Then we have}
\begin{enumerate} 
\item[{\bf 1.}] {\it the family $\{ g\0{t}\}$ is $\epsilon'$-slow, where
$\epsilon'=a'\epsilon$, with $a'=a'(n,\xi)$ }
\item[{\bf 2.}] {\it if $\epsilon<\frac{1}{\,3^22^{2}\,n!e^{2(1+\xi)}}$, then the family $\{ g\0{t} \}$ is $2$-bounded.}
\end{enumerate}
\noindent {\bf Remark.} We can take $a'=3e^{2(1+\xi)}\,(n+6\,c\0{14}(2))\,n^2\big(n!\,2^{n}  \big)^{3/2}$, where $c\0{14}(2)=
\frac{3}{2}n^{3/2}n!\,2^{n+2}$.\vspace{.1in}

\noindent {\bf Proof.} Recall that $|.|$ denotes the $C^2$ norm. The metric $g=e^{2t}g\0{t}+dt^2$ is $\epsilon$-close to hyperbolic means 
$$\big|e^{2t}\,(\,g\0{t}-\sigma\0{\R^n}\,) \big|\,=\,\big|(e^{2t}g\0{t}+dt^2)-(e^{2t}\sigma\0{\R^n}+dt^2)\big|\, =\, \big|g-\sigma\big|\, <\,\epsilon$$
\noindent And since $e^{2t}\geq e^{-2(1+\xi)}$, $t\geq-(1+\xi)$,
we get that $|g\0{t}-\sigma\0{\R^n}|_{C^0}=|\frac{1}{e^{2t}}(g-\sigma)|_{C^0}\leq e^{2(1+\xi)}\,\epsilon$. Using this same type of calculation
(for the derivatives of $g\0{t}-\sigma\0{\R^n}$)
together with the chain rule we obtain the following estimate 

\begin{equation} \big|\,g\0{t}-\sigma\0{\R^n}\, \big|\,\leq \,17 e^{2(1+\xi)}\,\epsilon\,< \,3^2 2\, e^{2(1+\xi)}\,\epsilon
\tag{i}\end{equation}

\noindent And, since we are assuming $2^{2}3^2\,n! e^{2(1+\xi)}\,\epsilon <1$ we get
$$ \big|g\0{t}\big|\, <\, \big|g\0{t}-\sigma\0{\R^n}\big|\, +\, \big|\sigma\0{\R^n}\big|\, \leq \,  (3^22 \, e^{2(1+\xi)}\epsilon)\, +\, 1 \, < \,1\,+\,1\,=\,2 $$
\noindent To prove that $|det\,g\0{t}|\0{C^0}>1/2$ we use the following
fact\vspace{.1in}

\noindent {\it Let $B=\{b_{ij}\}$ be an $n\times n$ matrix with
$max\{|b_{ij}|\}\leq b$. Then $det\, B>1-\sum_{k=1}^{n}
\dbinom{n}{k}(n-k)!\,b^k$. Hence if $b\leq 1$ then
$det\, B>1-2\, n!\,b.$}\vspace{.1in}

\noindent It follows from (i) that the matrices $g\0{t}$ can be written as $1\!\!1+B$, with $|B|<(3^22\, e^{2(1+\xi)}\,\epsilon )$.
(Here $1\!\!1=\sigma\0{\R^n}$ and $B=g\0{t}-\sigma\0{\R^n}$.) This
together with the fact above and the hypothesis $2^{2}3^2\,n! e^{2(1+\xi)}\,\epsilon <1$ imply that $det \,g\0{t}> 1- 2n! \,3^2 2e^{2(1+\xi)}>1/2$. \,\, This proves Statement (2).\\

We now prove Statement (1). Write $\tau=g-\sigma$.  Hence $g\0{t}=\sigma\0{\R^n}+e^{-2t}\tau$, with $|\tau|<\epsilon$. 
Hence for $u, v\in \R^n$ we have
$$ g\0{t} (u, v)\, =\, \langle u,v\rangle + e^{-2t}\tau (u,v)
$$
\noindent where $ \langle u,v\rangle=\sigma\0{\R^n}(u,v)$. Therefore
$ g'\0{t} (u, v)\, =\,  e^{-2t} \tau' (u,v)\, -\, 2\,e^{-2t}\tau (u,v)$. But
$|\tau|<\epsilon$, so taking
$u=\p_i,\, v=\p_j$ we get
 $|(g\0{ij})'\0{t}|\0{C^0}\leq 3e^{2(1+\xi)}\epsilon$. Similar calculations yield
$|(g\0{ij})''\0{t}|\0{C^0}\leq 9e^{2(1+\xi)}\epsilon$ and $|\frac{\p}{\p x\0{l}}(g\0{ij})'\0{t}|\0{C^0}\leq3e^{2(1+\xi)}\epsilon$.
The Corollary now follows from 3.6 by taking $c=2$ and
$\epsilon_1=9e^{2(1+\xi)}\epsilon$, $\epsilon_2=3e^{2(1+\xi)}\epsilon$.
This completes the proof of Corollary 3.7.\vspace{.1in}

As we mentioned in Section 2 our definition of radially $\epsilon$-close to hyperbolic metrics is not perfect.
For instance hyperbolic space is not radially $\epsilon$-close to hyperbolic near the
(chosen) center. The next result tells us how far we have to be from the center to
get the $\epsilon$-close to hyperbolicity of hyperbolic space. It follows from Corollary 3.3
and the fact that we can take $\epsilon=0$ in this particular case.\vspace{.1in}

\noindent {\bf Corollary 3.8.} {\it Let $o\in \HH^{n+1}$. Then hyperbolic $(n+1)$-space
$\HH^{n+1}$ is radially $\epsilon$-close to hyperbolic (with respect to $o$) outside
$\B_a(\HH^{n+1})$, with charts o excess $\xi$, provided}
$$
C_1' \, e^{-a}\,\leq \, \epsilon
$$ 

\noindent {\it where $C_1'=C_1'(n,\xi)$.}\vspace{.1in}

\noindent We can take $C_1'=C_1(c\0{\bS^n},n,\xi)$, with $C_1$ is as in 3.3, and $c\0{\bS^n}$ is such that $\sigma\0{\bS^n}$ is
$c\0{\bS^n}$-bounded.
We can rephrase this result in the following way.\vspace{.1in}

\noindent  {\bf Corollary 3.9.} {\it Let $o\in \HH^{n+1}$. Then hyperbolic $(n+1)$-space
$\HH^{n+1}$ is radially $\epsilon$-close to hyperbolic (with respect to $o$) outside
$\B_a(\HH^{n+1})$, with charts o excess $\xi$, provided
$ a\geq {\sf a}( \epsilon,n+1,\xi)$,
where  ${\sf a}(\epsilon,n+1,\xi)=ln\Big(\frac{C_1'}{\epsilon}  \Big)$. Here $C_1'=C_1(c\0{\bS^n},n,\xi)$, with $C_1$ is as in 3.3.}\\

\noindent {\bf \large 4.  The Two Variable Warping Deformation.}

Fix an atlas $\cA\0{\bS^n}$ on $\bS^n$ as before (see Remarks 1 and 2 after 3.1).
All norms and boundedness constants will be taken with respect to
this atlas.  Let $c\0{\bS^n}$ be a fixed constant such that $\sigma\0{\bS^n}$ is $c\0{\bS^n}$-bounded. 

Let $g$ be a metric on the $n$-sphere $\bS^n$ and consider the warped metric
$h=sinh^2t\, g +dt^2$ on $\bS^n\times \R^+$.
Also let $c\0{g}$ be such that $g$ is
$c\0{g}$-bounded. Write $c=c\0{g}+c\0{\bS^n}$.

Let  $\rho:\R\ra[0,1]$ be as in the Introduction. It can be checked that we can find such a $\rho$ with the following properties: (i) $|\rho'(t)|< 3$, (ii) $|\rho''(t)|<12$, for all $t$. Recall that in the Introduction we defined $\rho\0{a,d}(t)=\rho(2\,\frac{t-a}{d})$, for $a,\, d>0$. Then
$|\rho'\0{a,d}(t)|< 6/d$ \, and \, $|\rho''\0{a,d}(t)|<48/d^2$, for all $t$.

Consider the family of metrics $\big\{\sigma\0{\bS^n}+s(g-\sigma\0{\bS^n})\big\}\0{s\in[0,1]}$ on $\bS^n$. Note that this family is 
$c$-bounded. Also, since $[0,1]$ and $\bS^n$ are
compact, we have that there is $\epsilon\0{g}$  such that this
family is $\epsilon\0{g}$-slow (also see 3.6). 
By Lemma 3.1, Remark 6 after 3.1, and the fact that 
$|\rho'\0{a,d}(t)|,\,|\rho''\0{a,d}(t)|<6/d$ (assuming $d\geq 8$) we  have that \vspace{.1in}

\noindent {\bf (4.1)}\,\,\,\, The family of metrics $\big\{ \sigma\0{\bS^n} +\rho\0{a,d}(t) (g-\sigma\0{\bS^n})\big\}\0{t\in\R^+}$ is $c$-bounded and  $\Big(\frac{12}{d}\,\epsilon\0{g}\Big)$-slow. \vspace{.1in}

\noindent In the Introduction we defined  $(g\0{a,d})\0{t}=\sigma\0{\bS^n} +\rho\0{a,d}(t) (g-\sigma\0{\bS^n})$ and $\cT_{_{a,d}}\, g\, =\, sinh^2\, t\,\,(g\0{a,d})\0{t}+dt^2$. Hence, by 4.1, Corollary 3.3 and the definition of $\cT_{_{a,d}}\, g$ we get the following Corollary. \vspace{.1in}

\noindent {\bf Corollary 4.2.} {\it Let $g$ on $\bS^n$ be $c\0{g}$-bounded. Let $a,\, d>0$ and $b>1+\xi\geq 1$. Then $\cT_{_{a,d}}\, g\,$ is $\epsilon$-close to hyperbolic outside $\B_{b}$
with charts of excess $\xi$, provided} $$\,C_1\,\Big(e^{-b}+\frac{12}{d}\,\epsilon\0{g}\Big)\leq\epsilon$$
\noindent {\it where $C_1=C_1(c,n,\xi)$ (see 3.3),
$c=c\0{\bS^n}+c\0{g}$, and $\epsilon\0{g}$ is as above. } \vspace{.1in}

The next lemma says that the constant $\epsilon\0{g}$ depends only on the constant
$c\0{g}$. \vspace{.1in}

\noindent {\bf Lemma 4.3.}  {\it If $g$ is $c\0{g}$-bounded then we can take}
\begin{center}$\epsilon\0{g}=\epsilon\0{g}(c\0{g},n)\,=\,A\,\big(\,n\,, \,c\0{g}\,+\,c\0{\bS^n} \big)$\end{center}
\noindent {\it Here the function $A(n,x)$ is given by
\begin{center}$A(n,x)\,=\, x\, \Big( a(n,x')+b(n,x')\Big)$\end{center}
\noindent where $a$ and $b$ are as in 3.6 (see the remark after the statement of 3.6), and $x'=\big[n!\,x^{n+1}\big]^n$.}

\noindent {\bf Proof.} Write $\sigma=\sigma\0{\bS^n}$ and $c'=c\0{\bS^n}$.
We want to apply Lemma 3.6 to the family
$\{g\0{s}\}\0{s\in [0,1]}$, where $g\0{s}=\sigma+s(g-\sigma)=
(1-s)\sigma+sg$.  First we need the following claim.\vspace{.1in}

\noindent {\bf Claim.} {\it The family $\{g\0{s}\}\0{s\in [0,1]}$ is
$c''$-bounded, where $c''=\big[n!\,(c+c')^{n+1}\big]^n$.}\vspace{.07in}

\noindent First note that (recall $c>1$)
\begin{center}$|g\0{s}|\,\leq\, |g|\,+\,|\sigma|\,<\, 
c+c'\,<\,c''$\end{center}
\noindent It remains to prove that $det\, g\0{s}>\frac{1}{c''}$.
We will use the following fact: \vspace{.1in}

\noindent {\it Fact: Let $A$ be symmetric positive definite. Then
$u^TAu>d\,|u|^2$, for all $u$, if and only if all eigenvalues of $A$ are $>d$.} \vspace{.1in}
  
From the proof of Lemma 3.4 we get that all eigenvalues of
$g$ are $>\frac{1}{n!c^{n}}$ and 
all eigenvalues of $\sigma$ are $>\frac{1}{n!(c')^{n}}$.
Hence all eigenvalues of either $g$ or $\sigma$ are $>e=\frac{1}{n!(c+c')^{n}}$. 
Therefore (using the fact above) we get $h(u,u)> e\, |u|^2$, for every $u\in\R^n$, where $h$ is either
$g$ or $\sigma$. Then
$$
g\0{s}(u,u)\,=\,(1-s)\,\sigma(u,u)\,+\,s\,g(u,u)\,>\,e\,|u|^2
$$
\noindent which implies, by the fact above, that all eigenvalues of $g\0{s}$
are $>e$. Thus $det \,g\0{s}>e^n=c''$. This proves the claim. \vspace{.1in}

To finish the proof of the lemma apply Lemma 3.6 the the family
$\{g\0{s}\}\0{s\in [0,1]}$, which, by the claim, is $c''$-bounded.
A simple calculation shows that in this case we can take $\epsilon\0{1}$
and  $\epsilon\0{2}$ in Lemma 3.6 satisfying
 $\epsilon\0{1}=\epsilon\0{2}=|\sigma|_{C^2}+|g|_{C^2}<c+c'$.
This proves the lemma. \vspace{.1in}

The next result is the slightly more general version of Theorem
1 in the Introduction.
\vspace{.1in}

\noindent {\bf Corollary 4.4.} {\it Let the metric $g$ on $\bS^n$ be $c\0{g}$-bounded. Let $a,\, d>0$, $b>1+\xi\geq 1$. Then \vspace{.1in}

\noindent (1) the metric $\cT_{_{a,d}}\, g\,$ is 
hyperbolic on $B_{a}$ 

\noindent (2) the metric $\cT_{_{a,d}}\, g\,$ is 
radially $\epsilon$-close to hyperbolic outside $B_{b}$,
with charts of excess $\xi$, \vspace{.1in}
 
\noindent provided} 
\begin{center}$\,C_2\,\Big(e^{-b}+\frac{1}{d}\,\Big)\leq\epsilon$
\end{center}
\noindent {\it where $C_2=C_2(c\0{g},n,\xi)$. } \vspace{.1in}

\noindent {\bf Proof.} First note that by Corollary 4.3
we have $\epsilon\0{g}=\epsilon\0{g}(c\0{g},n)$. Also note that the maximum
value of $\frac{e^{-b}+\frac{12\epsilon\0{g}}{d}}{e^{-b}+\frac{1}{d}}$, $d,\, b>0$,
is less than $1+12\epsilon\0{g}$. Hence if we take $C_2=C_2(c,n,\xi)=(1+12\epsilon\0{g}(c,n))C_1(c,n,\xi)$ we get
$C_2\, (e^{-b}+\frac{1}{d})\geq C_1\, (e^{-b}+\frac{12\epsilon\0{g}}{d})$.
The result now follows from Corollary 4.2. This proves the corollary.\vspace{.1in}

Note that Theorem 1 is obtained from Corollary 4.4 by taking
the excess $\xi=0$. Let
$B_{a}=B_a(0)$ be the ball of radius $a$ centered at $0$. We say
that a metric $h$ on $\R^{n+1}$ is $(B_a,\epsilon)$-{\it close to
hyperbolic, with charts of excess $\xi$}, if \vspace{.1in}

\begin{enumerate}
\item[  (1)]  On $B_{a}-\{0\}=\bS^n\times (0,a)$
we have $h=sinh^2(t)\sigma\0{\bS^n}+dt^2$. Hence $h$
is hyperbolic on $B_a$.
\item[(2)]  the metric $h$ is 
radially $\epsilon$-close to hyperbolic outside $B_{a-1-\xi}$,
with charts of excess $\xi$.
\end{enumerate}\vspace{.1in}
  
\noindent {\bf Remarks.}

\noindent {\bf 1.} We have dropped the word ``radially" 
to simplify the notation. But it does appear in condition (2),
where now ``radially" refers to the center on $B_a$.

\noindent{\bf 2.} We will always assume $a>{\sf a}+1+\xi$,
where {\sf a} is as in 3.9. Therefore conditions (1), (2) 
and 3.9 imply
a stronger version of (2):\vspace{.1in}

 (2') the metric $h$ is 
radially $\epsilon$-close to hyperbolic outside $B_{\sf a}$,
with charts of excess $\xi$.\vspace{.1in}

\noindent This is the reason why we demanded radius $a-1-\xi$
in (2), instead of just $a$. \vspace{.1in}

With this new notation, Corollary 4.4 can be restated in the following way (taking
$b=a-1-\xi$ in 4.4): \vspace{.1in}

\noindent {\bf Corollary 4.5.} {\it Let the metric $g$ on $\bS^n$ be $c$-bounded. Then the metric $\cT_{_{a,d}}\, g\,$ is  $(B_a,\epsilon)$-close to hyperbolic, with charts of excess $\xi$, provided} 
\begin{center}$\,C_2\,\Big(e^{-a}+\frac{1}{d}\,\Big)\leq\epsilon$
\end{center}
\noindent {\it where $C_2$ is as in 4.4. } \vspace{.1in}

Note that Corollary in the Introduction is obtained from 4.5
by taking $\xi=0$.
The next result is the slightly more general version of Theorem
2 in the Introduction. It follows directly from 4.5 by taking\vspace{.1in}

$$\begin{array}{rcl}  a(c,\epsilon,n,\xi)&=&ln(\frac{2C_2}{\epsilon})+{\sf a}(\epsilon, n+1,\xi)\\ \\d(c,\epsilon,n,\xi)&=&\frac{2C_2}{\epsilon}\end{array}$$\vspace{.1in}

\noindent {\bf Corollary 4.6.} {\it Let the metric $g$ on $\bS^n$ be $c$-bounded and $ \epsilon>0$ and $\xi\geq 0$. Then the metric $\cT_{_{a,d}}\, g\,$ is 
 $(B_a,\epsilon)$-close to hyperbolic,
with charts of excess $\xi$, provided we take $a$ and $d$ large enough. Explicitly we have to take} 
\begin{center}$a\,>\, a(c,\epsilon,n,\xi)\hspace{.4in} and\hspace{.4in}  d\,>\, d
(c,\epsilon,n,\xi)$\end{center}\vspace{.1in}

Note that Theorem 2 in the Introduction is obtained from Corollary 4.6 by taking the excess $\xi=0$. 

\vspace{.1in}

{\footnotesize

Pedro Ontaneda

SUNY, Binghamton, N.Y., 13902, U.S.A.}


\begin{thebibliography}{99}






\bibitem{AF1} C. S. Aravinda and F.T. Farrell, 
{\em Rank 1 aspherical manifolds which do not support any nonpositively curved metric},
Comm. in Analysis and Geometry {\bf 2} (1994) 65-78.




\bibitem{AF2} C. S. Aravinda and F.T. Farrell, 
{\em  Exotic negatively curved structures on Cayley hyperbolic manifolds},
Jour. Diff. Geom. {\bf 63} (2003) 41-62.





\bibitem{AF3} C. S. Aravinda and F.T. Farrell, 
{\em Exotic structures and quaternionic hyperbolic manifolds}.
Algebraic groups and Arithmetic, (Eds: S. G. Dani and G. Prasad)
Narosa Publishing House, 2004.


\bibitem{A} S. Ardanza-Trevillano, {\em Exotic smooth structures on negatively
curved manifolds that are not of the homotopy type of a locally symmetric space}.
PhD. Thesis, SUNY Binghamton (2000).




\bibitem{BisOn} R.L. Bishop and B. O'Neill, {\em Manifolds of negative 
curvature}, Trans. Amer. Math. Soc. {\bf145} (1969) 1-49. 







   











\bibitem{ChD} R. M. Charney and M. W. Davis, 
{\em Strict hyperbolization},
Topology {\bf 34} (1995), 329-350. 

































\bibitem{FJ1} F.T. Farrell and L.E. Jones, {\em Negatively curved manifolds with exotic smooth structures}, J. Amer. Math. Soc. 
{\bf 2} (1989) 899-908.


\bibitem{FJ1.5}  F.T. Farrell and L.E. Jones, {\em Nonuniform hyperbolic lattices and exotic smooth structures}, J. Differential Geom., {\bf 38} (1993), 235-261.




\bibitem{FJ2} F.T. Farrell and L.E. Jones, {\em  Complex hyperbolic manifolds and exotic smooth structures}, Invent. Math., {\bf 117} (1994), 57-74.




\bibitem{FJO1} F.T. Farrell, L.E. Jones and P. Ontaneda, {\em Hyperbolic manifolds with negatively curved exotic 
triangulations in dimension larger than five}. Jour. Diff. Geom. {\bf 48} (1998) 319-322.

\bibitem{FJO2} F.T. Farrell, L.E. Jones and P. Ontaneda, {\em Examples of non-homeomorphic
harmonic maps between negatively curved manifolds}, Bull. London Math. Soc. {\bf 30} 
(1998) 295-296.

\bibitem{FOR} F.T. Farrell, P. Ontaneda and M.S. Raghunathan, {\em Non-univalent harmonic maps homotopic
to diffeomorphisms}, Jour. Diff. Geom. {\bf 54} (2000) 227-253.

\bibitem{FO1} F.T. Farrell and P. Ontaneda, {\em Cellular harmonic maps which are not diffeomorphisms}.
Inventiones Mathematicae {\bf 158 (1)} 497-513, (2004).


\bibitem{FO1.5} F.T. Farrell and P. Ontaneda, {\em Branched cover of hyperbolic manifolds and harmonic maps}, Comm. in Analysis and Geometry, {\bf 14}(2) (2006), 249�268.

\bibitem{FO2} F.T. Farrell and P. Ontaneda, {\em A caveat on the convergence of the Ricci flow
for negatively curved manifolds}, Asian Journal of Mathematics, {\bf  9 (3)}, 401-406, (2005).



\bibitem{FO3} F.T. Farrell and P. Ontaneda, {\em The Teichm\"{u}ller Space of Pinched Negatively Curved Metrics on a Hyperbolic Manifold is not Contractible}, Annals of Mathematics (2) {\bf 170} (2009), 45-65.






\bibitem{FO4} F.T. Farrell and P. Ontaneda, {\em Teichm\"uller spaces and bundles with negatively curved fibers}, GAFA {\bf 20} (2010), 1397-1430.


\bibitem{FO5} F.T. Farrell and P. Ontaneda, 
{\em The Moduli Space of Negatively Curved Metrics of a Hyperbolic Manifold}, Journal of Topology {\bf 3} (2010), 561-577. 












































\bibitem{O} P. Ontaneda, {\em Hyperbolic manifolds with negatively curved exotic triangulations in dimension six}, 
J. Diff. Geom. {\bf 40} (1994), 7-22. 

\bibitem{O1} P. Ontaneda, {\em  Riemannian Hyperbolization}.  	ArXiv:
1406.1730.

\bibitem{O2} P. Ontaneda, {\em  Hyperbolic Extensions and Metrics $\epsilon$-Close to Hyperbolic}.  ArXiv: 1406.1740.


\bibitem{O3} P. Ontaneda, {\em  Deforming an $\epsilon$-Close to Hyperbolic metric to a warp metric}.  ArXiv:1406.1741.


\bibitem{O4} P. Ontaneda, {\em Deforming an $\epsilon$-Close to Hyperbolic Metric to a Hyperbolic Metric}. ArXiv: 1406.1743.











































\end{thebibliography}
\end{document}